\documentclass[a4paper,12pt]{article}
\usepackage{amssymb,graphicx}

\newtheorem{theorem}{Theorem}

\begin{document}

\author{
Dumitru B\u{a}leanu\\
\small{\c{C}ankaya University,}\\
\small{Department of Mathematics {\&} Computer Science,}\\
\small{\"{O}gretmenler Cad. 14 06530, Balgat -- Ankara, Turkey,}\\
\small{{\&} Institute of Space Sciences, M\u{a}gurele -- Bucure\c{s}ti, Romania}\\
\small{e-mail address: dumitru@cankaya.edu.tr}\\
Octavian G. Mustafa\\
\small{University of Craiova, }\\
\small{Faculty of Exact Sciences,}\\
\small{A{.}~I{.} Cuza 13, 200534 Craiova, Romania}\\
\small{e-mail address: octaviangenghiz@yahoo.com}\\
and\\
Donal O'Regan\\
\small{National University of Ireland,}\\
\small{School of Mathematics, Statistics and Applied Mathematics,}\\
\small{Galway, Ireland}\\
\small{e-mail address: donal.oregan@nuigalway.ie}
}

\title{A Kamenev-type oscillation result for a linear\\
 $(1+\alpha)$--order fractional differential equation}
\date{}
\maketitle

\noindent{\bf Abstract} We investigate the eventual sign changing for the solutions of the linear equation $\left(x^{(\alpha)}\right)^{\prime}+q(t)x=0$, $t\geq0$, when the functional coefficient $q$ satisfies the Kamenev-type restriction $\limsup\limits_{t\rightarrow+\infty}\frac{1}{t^{\varepsilon}}\int_{t_0}^{t}(t-s)^{\varepsilon}q(s)ds=+\infty$ for some $\varepsilon>2$, $t_{0}>0$. The operator $x^{(\alpha)}$ is the Caputo differential operator and $\alpha\in(0,1)$.

\section{Introduction}
The oscillation of solutions for ordinary differential equations is
an important topic in applied mathematics. We note that the KBM
(Krylov-Bogoliubov-Mitropolsky) averaging technique and the theory of
adiabatic invariants were applied successfully to problems in
celestial mechanics \cite[pp{.} 41, 195]{boccaletti} that can be connected with the oscillation theory.

In the particular case of the second order linear differential equation
\begin{eqnarray}
x^{\prime\prime}+q(t)x=0,\quad t\geq0,\label{ode}
\end{eqnarray}
where the functional coefficient
$q:[0,+\infty)\rightarrow\mathbb{R}$ is continuous, I{.}V{.} Kamenev
\cite{kamenev} proved in 1978 that oscillations occur when
\begin{eqnarray}
\limsup\limits_{t\rightarrow+\infty}\frac{1}{t^{\varepsilon}}\int_{t_0}^{t}(t-s)^{\varepsilon}q(s)ds=+\infty\label{kam}
\end{eqnarray}
for some $\varepsilon>1$ and $t_{0}>0$. This result replaces the
classical Wintner-Hartman averaging quantity
$\limsup\limits_{t\rightarrow+\infty}\frac{1}{t}\int_{t_0}^{t}\int_{t_0}^{s}q(\tau)d\tau
ds$ from the various statements of oscillation criteria regarding (\ref{ode}) with the left-hand part of (\ref{kam}). In the original paper the number $\varepsilon\geq2$ was an
integer, but J{.}S{.}W{.} Wong \cite[pp{.} 418--419]{wong} noticed
that it could be recast with any real number greater than $1$.

The aim of this paper is to present a Kamenev type theorem in the
framework of fractional differential equations. To the best of our
know\-ledge, such a result has not been established for any
generalized differential equation.

Differential equations of non-integer order, al\-so called \textit{frac\-tio\-nals} {(}FDE{'} s{)},  a\-rise naturally in models in en\-gi\-ne\-e\-ring, phy\-sics or che\-mis\-try and we refer the reader to
\cite{baleanuetal,diethelm,trujillo,miller,podlubny}.

Consider a function $h\in C^{1}(I,\mathbb{R})\cap
C(\overline{I},\mathbb{R})$ with
$\lim\limits_{t\searrow0}[t^{1-\alpha}h^{\prime}(t)]\in\mathbb{R}$
for some $\alpha\in(0,1)$, where $I=(0,+\infty)$. The \textit{Caputo
derivative} of order $\alpha$ of $h$ is defined as
\begin{eqnarray*}
h^{(\alpha)}(t)=\frac{1}{\Gamma(1-\alpha)}\cdot\int_{0}^{t}\frac{h^{\prime}(s)}{(t-s)^{\alpha}}ds,\quad t\in I,
\end{eqnarray*}
where $\Gamma$ is the Euler function Gamma, cf{.} \cite[p{.}
79]{podlubny}. Note if we let the function $h^{\prime}$ be
absolutely continuous \cite[Chapter 7]{rudin} then the (usual)
derivative of $h^{(\alpha)}$ will exist almost everywhere with
respect to the Lebesgue measure $m$ on $\mathbb{R}$, see \cite[p{.}
35, Lemma 2{.}2]{samko}. Further, we have that
\begin{eqnarray}
h(t)=h(0)+\frac{1}{\Gamma(\alpha)}\int_{0}^{t}\frac{h^{(\alpha)}(s)}{(t-s)^{1-\alpha}}ds,\quad t\in I,\label{integro}
\end{eqnarray}
provided that $h^{(\alpha)}$ is in $L^{\infty}(m)$.

The FDE we investigate in this paper is
\begin{eqnarray}
\left(x^{(\alpha)}\right)^{\prime}(t)+q(t)x(t)=0,&t\in I,\label{iniprob}
\end{eqnarray}
where the continuous functional coefficient
$q:\overline{I}\rightarrow\mathbb{R}$ satisfies the Kamenev
condition  (\ref{kam}) for some $\varepsilon>2$. The asymptotic
behavior of solutions to (\ref{iniprob}) was discussed in \cite{bm2}
and the authors showed that if
\begin{eqnarray}
\int_{0}^{+\infty}t^{1+\alpha}\vert q(t)\vert dt<+\infty\quad\mbox{and}\quad\int_{0}^{+\infty}t^{\alpha}\vert q(t)\vert dt<\Gamma(1+\alpha),\label{asympt_infty}
\end{eqnarray}
then, for every $c_{1},\thinspace c_{2}\in\mathbb{R}$, the equation (\ref{iniprob}) has a solution $x$ with the asymptotic expression
\begin{eqnarray}
x(t)=c_{1}+c_{2}\cdot t^{\alpha}+o(1)\quad\mbox{when }t\rightarrow+\infty.\label{asympt1}
\end{eqnarray}

Finally we mention a recent contribution \cite{grace} which concerns
oscillation of perturbed FDE's with power-like nonlinearities. The
proofs there rely exclusively on the averaging of the perturbation thus being completely different from the method in our investigation.

\section{Statement of our result and a comment}

Throughout this note, by \textit{a solution} to the
$(1+\alpha)$--order FDE (\ref{iniprob}) we mean any function $x\in
C^{1}(\overline{I},\mathbb{R})$ that verifies (\ref{iniprob}) in
$I$. Such a solution $x$ \textit{oscillates} if there exists an
increasing, unbounded from above, sequence $(t_{n})_{n\geq1}\subset
I$ such that
\begin{eqnarray*}
x(t_{2n-1})<0\quad\mbox{and}\quad x(t_{2n})>0,\quad n\geq1.
\end{eqnarray*}

\begin{theorem}\label{kam_th}
Any solution $x$ of (\ref{iniprob}), (\ref{kam}) either oscillates or satisfies the ine\-quality
\begin{eqnarray}
\liminf\limits_{t\rightarrow+\infty}\left\{x^{(\alpha)}(t)\cdot\left[x^{\prime}(t)-x^{(\alpha)}(t)\right]\right\}\leq0.\label{sign_chang_ineq}
\end{eqnarray}
More precisely, in the situation (\ref{sign_chang_ineq}), there is an increasing, unbounded from above, sequence $(T_{n})_{n\geq1}\allowbreak\subset I$ such that
\begin{eqnarray}
x^{(\alpha)}(T_{n})\cdot\left[x^{\prime}(T_{n})-x^{(\alpha)}(T_{n})\right]<0,\quad n\geq1.\label{sign_chang_ineq_1}
\end{eqnarray}
\end{theorem}

At first glance, the conclusion of our result is rather
disappointing given the fact that we are  not able to insulate
oscillations from other types of asymptotic behavior. However, let
us recall the classical Fite oscillation criterion \cite{fite} where
the conclusion is again formulated as a list of multiple outcomes.

To establish that \textit{the possibility (\ref{sign_chang_ineq}) cannot be removed from the statement}, let us consider the case when $q(t)>0$ for all $t\geq0$. Obviously, from the in\-e\-qua\-li\-ties
\begin{eqnarray*}
\frac{1}{t^{\varepsilon}}\int_{t_0}^{t}(t-s)^{\varepsilon}q(s)ds\leq\frac{1}{t^{\varepsilon}}\int_{t_0}^{t}t^{\varepsilon}q(s)ds\leq\int_{t_0}^{+\infty}q(s)ds,\quad t_{0}>0,
\end{eqnarray*}
we get that $\int_{0}^{+\infty}q(t)dt=+\infty$.

Assume now that $x$ is a non-oscillatory solution to
(\ref{iniprob}), which implies, without loss of generality, that we
can take $x(t)>0$ for every $t\geq T>0$. Since
$\left(x^{(\alpha)}\right)^{\prime}(t)<0$ in $[T,+\infty)$, that is,
the function $x^{(\alpha)}$ is eventually decreasing, there exists
$L\in[-\infty,+\infty)$ such that
$\lim\limits_{t\rightarrow+\infty}x^{(\alpha)}(t)=L$.

Suppose further, for the sake of contradiction, that (\ref{sign_chang_ineq}) does not hold either, i.e.
\begin{eqnarray}
x^{(\alpha)}(t)\cdot\left[x^{\prime}(t)-x^{(\alpha)}(t)\right]\geq0\quad\mbox{for all}\quad t\geq T.\label{contra}
\end{eqnarray}

Consider first the case $L<0$. Now since $x^{(\alpha)}$ becomes
eventually negative valued then the function
$x^{\prime}-x^{(\alpha)}$ becomes eventually non-positive valued.
Thus there exists a $T_{1}\geq T$ large enough so that
\begin{eqnarray*}
x^{\prime}(t) \leq x^{(\alpha)}(t)<\frac{L}{2},\quad t\geq T_{1}.
\end{eqnarray*}

An integration with respect to the variable $t$ leads to $x(t)\leq x(T_{1})+\frac{L}{2}\cdot(t-T_{1})$
and so $\lim\limits_{t\rightarrow+\infty}x(t)=-\infty$, which contradicts the eventual positivity of $x(t)$.

Consider next the case $L>0$. Now, $x^{\prime}(t)\geq
x^{(\alpha)}(t)>\frac{L}{2}$ for any $t\geq T_{2}\geq T$ large
enough. We get that $\lim\limits_{t\rightarrow+\infty}x(t)=+\infty$
and also, as a by-product, $\int_{0}^{+\infty}q(t)x(t)dt=+\infty$.

However, since
\begin{eqnarray*}
x^{(\alpha)}(t)=x^{(\alpha)}(T_{2})-\int_{T_{2}}^{t}q(s)x(s)ds,\quad t\geq T_{2},
\end{eqnarray*}
we deduce that $\lim\limits_{t\rightarrow+\infty}x^{(\alpha)}(t)=-\infty$ which, again, contradicts our hypotheses.

Finally consider the case $L=0$. Since the function $x^{(\alpha)}$
is eventually decreasing, we have $x^{(\alpha)}(t)>L=0$ for all
$t\geq T$. Similarly, $x^{\prime}(t)\geq x^{(\alpha)}(t)>0$
throughout $[T,+\infty)$. This yields $x(t)\geq x(T)>0$ for all
$t\geq T$ and so
\begin{eqnarray*}
x^{(\alpha)}(t)&=&x^{(\alpha)}(T)-\int_{T}^{t}q(s)x(s)ds\leq x^{(\alpha)}(T)-x(T)\cdot\int_{T}^{t}q(s)ds\\
&\rightarrow&-\infty\quad\mbox{when }t\rightarrow+\infty,
\end{eqnarray*}
a contradiction.

What kind of functions verify (\ref{sign_chang_ineq})? An elementary
example --- though not from $C^{1}(\overline{I},\mathbb{R})$ --- is
$x$ with $x(t)=t^{\beta}$, $t\geq0$, for some $\beta\in(0,\alpha)$.
Here,
$x^{(\alpha)}(t)=\frac{\Gamma(1+\beta)}{\Gamma(1+\beta-\alpha)}\cdot
t^{\beta-\alpha}$, $t>0$. The coefficient $q(t)$ of (\ref{iniprob})
reads now as $q(t)=C(\alpha,\beta)\cdot t^{-1-\alpha}$, where
$C(\alpha,\beta)=\frac{(\alpha-\beta)\Gamma(1+\beta)}{\Gamma(1+\beta-\alpha)}$,
and, unfortunately, does not satisfy the condition (\ref{kam}),
since
\begin{eqnarray*}
\frac{1}{t^{\varepsilon}}\int_{t_0}^{t}(t-s)^{\varepsilon}q(s)ds&=&\left(\frac{t-t_{0}}{t}\right)^{\varepsilon}\int_{t_0}^{+\infty}q(\tau)d\tau\\
&-&\frac{\varepsilon}{t^{\varepsilon}}\int_{t_0}^{t}\left(\int_{s}^{+\infty}q(\tau)d\tau\right)\frac{ds}{(t-s)^{1-\varepsilon}}\\
&\leq&\int_{t_0}^{+\infty}q(\tau)d\tau<+\infty.
\end{eqnarray*}

Even though the functional coefficient $q(t)=C(\alpha,\beta)\cdot
t^{-1-\alpha}$ does not satisfy either of the restrictions in
(\ref{asympt_infty}), it seems to us that there is no easy way to
determine closed form solutions of (\ref{iniprob}), (\ref{kam}) that
will obey (\ref{sign_chang_ineq_1}). On the other hand, notice that,
for any positive constant $A$, the functional coefficient $q(t)=A$
verifies the Kamenev condition (\ref{kam}). The formula of the
Laplace transform for Caputo derivatives \cite[p{.} 106, Eq{.}
(2{.}253); p{.} 21, Eq{.} (1{.}80)]{podlubny} leads us to the
oscillatory solution $x(t)=E_{1+\alpha}(-At^{1+\alpha})$, $t\geq0$,
of (\ref{iniprob}), where $E_{\gamma}$ denotes the Mittag-Leffler
function \cite[p{.} 16]{podlubny}.

\section{Proof of Theorem \ref{kam_th}}

Assume that the solution $x$ of (\ref{iniprob}) does not possess any zeros in $[T,+\infty)$ for some $T\geq0$ large enough. Suppose also, for the sake of contradiction, that (\ref{contra}) holds true.

We introduce the quantity $w(t)=\frac{x^{(\alpha)}(t)}{x(t)}$, where
$t\geq T$. Now we have
\begin{eqnarray*}
w^{\prime}(t)=-q(t)-\frac{x^{(\alpha)}(t)\cdot x^{\prime}(t)}{[x(t)]^{2}},\quad t\geq T,
\end{eqnarray*}
and so we have the typical Riccati inequality
\begin{eqnarray}
w^{\prime}(t)+[w(t)]^{2}+q(t)=\frac{-x^{(\alpha)}(t)\cdot x^{\prime}(t)+\left[x^{(\alpha)}(t)\right]^{2}}{[x(t)]^{2}}\leq0,\label{ricc}
\end{eqnarray}
valid for every $t\geq T$.

Further, we deduce that
\begin{eqnarray}
\int_{T}^{t}(t-s)^{\varepsilon}q(s)ds&\leq&-\int_{T}^{t}(t-s)^{\varepsilon}\left\{w^{\prime}(s)+[w(s)]^{2}\right\}ds\nonumber\\
&=&w(T)\cdot(t-T)^{\varepsilon}-\varepsilon\int_{T}^{t}w(s)(t-s)^{\varepsilon-1}ds\nonumber\\
&-&\int_{T}^{t}(t-s)^{\varepsilon}[w(s)]^{2}ds\nonumber\\
&\leq&\vert w(T)\vert\cdot(t-T)^{\varepsilon}+\varepsilon\int_{T}^{t}\vert w(s)\vert(t-s)^{\varepsilon-1}ds\nonumber\\
&-&\int_{T}^{t}(t-s)^{\varepsilon}[w(s)]^{2}ds.\label{step1}
\end{eqnarray}

Notice as well that
\begin{eqnarray}
&&(t-s)^{\varepsilon}[w(s)]^{2}-\varepsilon\vert w(s)\vert(t-s)^{\varepsilon-1}\nonumber\\
&&=\left[(t-s)^{\frac{\varepsilon}{2}}\vert w(s)\vert -\frac{\varepsilon}{2}(t-s)^{\frac{\varepsilon}{2}-1}\right]^{2}-\frac{\varepsilon^{2}}{4}(t-s)^{\varepsilon-2},\label{step2}
\end{eqnarray}
where $t\geq s\geq T$.

By combining (\ref{step1}), (\ref{step2}), we are able to eliminate the quantity $w(s)$ from the estimate of $q$, namely
\begin{eqnarray*}
\int_{T}^{t}(t-s)^{\varepsilon}q(s)ds\leq\vert w(T)\vert(t-T)^{\varepsilon}+\frac{\varepsilon^{2}}{4}\cdot\frac{(t-T)^{\varepsilon-1}}{\varepsilon-1}.
\end{eqnarray*}

Thus,
\begin{eqnarray*}
\frac{1}{t^{\varepsilon}}\int_{T}^{t}(t-s)^{\varepsilon}q(s)ds\leq\vert w(T)\vert+\frac{\varepsilon^{2}}{4(\varepsilon-1)}\cdot\frac{1}{t},\quad t\geq T.
\end{eqnarray*}
This estimate, obviously, contradicts the Kamenev condition (\ref{kam}).

The proof is complete.

\bigskip

Let us return to (\ref{sign_chang_ineq}) for other comments.

Firstly, suppose that equation (\ref{iniprob}) has a solution $x\in
C^{2}(\overline{I},\mathbb{R})$ such that
\begin{eqnarray}
x^{\prime\prime}(t)\leq0,\thinspace t\geq T>0,\quad\mbox{and}\quad\lim\limits_{t\rightarrow+\infty}x^{\prime}(t)=L\in I.\label{asympt2}
\end{eqnarray}
Without loss of generality, we may assume that
\begin{eqnarray*}
x^{\prime\prime}(t)\leq0,\quad\frac{3L}{2}\geq x^{\prime}(t)\geq\frac{L}{2},\quad t\geq T.
\end{eqnarray*}

For any $t\geq 2T$, we deduce that
\begin{eqnarray*}
&&\frac{1}{\Gamma(1-\alpha)}\int_{0}^{t}\frac{x^{\prime}(s)}{(t-s)^{\alpha}}ds=\frac{1}{\Gamma(1-\alpha)}\left(\int_{0}^{T}+\int_{T}^{t}\right)\frac{x^{\prime}(s)}{(t-s)^{\alpha}}ds\\
&&\geq\frac{1}{\Gamma(1-\alpha)}\left[\int_{T}^{t}\frac{x^{\prime}(s)}{(t-s)^{\alpha}}ds-\frac{1}{(t-T)^{\alpha}}\int_{0}^{T}\vert x^{\prime}(s)\vert ds\right]\\
&&\geq\frac{1}{\Gamma(1-\alpha)}\left[x^{\prime}(t)\cdot\int_{T}^{t}\frac{1}{(t-s)^{\alpha}}ds-\frac{c(x,T)}{(t-T)^{\alpha}}\right]\\
&&\geq\frac{1}{\Gamma(1-\alpha)}\left[\frac{L}{2}\cdot\frac{(t-T)^{1-\alpha}}{1-\alpha}-\frac{c(x,T)}{T^{\alpha}}\right],
\end{eqnarray*}
where $c(x,T)=\int_{0}^{T}\vert x^{\prime}(s)\vert ds=x(T)-x(0)$.

Thus,
\begin{eqnarray*}
x^{(\alpha)}(t)\geq2L>\frac{3L}{2}\geq x^{\prime}(t)\geq\frac{L}{2}>0,\quad t\geq T_{3}\geq 2T,
\end{eqnarray*}
for some $T_{3}=T_{3}(\alpha,x,T)$ large enough. As a by-product,
\begin{eqnarray*}
x^{(\alpha)}(t)\cdot[x^{\prime}(t)-x^{(\alpha)}(t)]<0,\quad t\geq T_{3}.
\end{eqnarray*}

In conclusion, \textit{if the Kamenev condition (\ref{kam}) would allow the existence of solutions to (\ref{iniprob}) verifying (\ref{asympt2}) then these solutions are candidates for the estimate (\ref{sign_chang_ineq})}.

To make a connection with (\ref{asympt1}), notice that the solutions from (\ref{asympt2}) have the asymptotic expression
\begin{eqnarray*}
x(t)=c_{3}\cdot t+o(t)\quad\mbox{when }t\rightarrow+\infty,
\end{eqnarray*}
where $c_{3}\in\mathbb{R}$. Such solutions, usually called
\textit{asymptotically linear}, are of in\-te\-rest in the theory of
fractional differential equations, see \cite{bm1}.

An open question is whether such solutions exist. However, in the
case of equation (\ref{ode}), if the functional coefficient $q$
satisfies the restriction
\begin{eqnarray*}
\int_{0}^{+\infty}t\cdot\max\{q(t),0\}dt=+\infty,
\end{eqnarray*}
then for all solutions $x$ we get
$\lim\limits_{t\rightarrow+\infty}x^{\prime}(t)=0$, meaning they
cannot verify (\ref{asympt2}). Also, if
\begin{eqnarray*}
\int_{0}^{+\infty}-\min\{q(t),0\}dt<+\infty,
\end{eqnarray*}
then for all solutions $x$ we obtain
$x^{\prime}(t)=O\left(t^{-\frac{1}{2}}\right)$ when
$t\rightarrow+\infty$; see \cite{boas,wintner}.

Even though we have focused here on FDE's, the principle in Theorem
\ref{kam_th} may be applied to various nonlinear differential
equations. For example, assume that we replace the Caputo derivative
$x^{(\alpha)}$ with ${\cal D}x=\frac{x^{\prime}}{\sqrt{1+(
x^{\prime})^{2}}}$. It is clear that
\begin{eqnarray*}
\left({\cal D}x\right)(t)\cdot\left[x^{\prime}(t)-\left({\cal D}x\right)(t)\right]&=&\frac{[x^{\prime}(t)]^{2}}{1+[x^{\prime}(t)]^{2}}\cdot\left\{\sqrt{1+[x^{\prime}(t)]^{2}}-1\right\}\\
&\geq&0,\quad t\in I.
\end{eqnarray*}

According to Theorem \ref{kam_th}, whenever the functional coefficient $q$ obeys the Kamenev restriction (\ref{kam}), the non-trivial solutions of the differential equation
\begin{eqnarray*}
\left({\cal D}x\right)^{\prime}+q(t)x=0,\quad t>0,
\end{eqnarray*}
oscillate. This conclusion complements the foundational work in \cite[Sect{.} 4{.}3]{pucci_serrin}.

\textbf{Acknowledgment.} The work of D{.} B{.} and O{.} G{.} M{.} has been supported by a grant of (to be completed by Prof. B\u{a}leanu).

\end{document}